\documentclass[a4,10pt]{amsart}
\usepackage{amssymb,amsmath,graphicx,tikz}

\def\H{\mathbb{H}}
\def\D{\mathbb{D}}

\def\C{\mathbb{C}}

\def\Z{\mathbb{Z}}

\def\Re{{\rm Re} \,}
\def\Im{{\rm Im} \,}

\def\2int{\mathop{\int\int}}

\def\kzu1#1{\buildrel #1:1 \over \longrightarrow}

\def\ds{\displaystyle}
\def\ts{\textstyle}

\def\j{\mathfrak{j}}

\def\0{\mathbf{0}}
\def\1{\mathbf{1}}

\def\be#1{\begin{equation}\label{#1}}
\def\ee{\end{equation}}

\begin{document}

\title{Conformal Mapping of Circular Triangles \& Classical Function Theory}
\author{Norbert Steinmetz}
\maketitle

{\small
\begin{abstract} In the present paper it is discussed to which extend conformal mappings of special circular triangles, which are thoroughly examined by M.\ Bonk~\cite{bonk1,bonk2} using analytic
and classical special functions methods, can be defined and investigated geometrically.
This is in analogy to the elliptic modular function, which admits an analytic and also a geometric definition. \end{abstract}

\bigskip
\begin{tabular}{rl}
{\sc Keywords.}& Circular triangle, conformal mapping, Weierstrass P- and zeta-function,\cr
& $\j$-function, Schwarzian derivative\cr
{\sc 2020 MSC.}& 30C20, 33E05, 34M04 
\end{tabular}}

\section{Introduction}
\begin{figure}[!ht]
\begin{center}
\begin{tabular}{cc}
\begin{tikzpicture}
[scale=0.6];
\draw[line width=0.5,white] (0,0)--(0,-2.5);\draw[line width=0.5,gray]  (0,0) circle (2);\draw[line width=0.5,gray] (0,5)--(0,0);\draw[line width=0.5,gray] (1,5)--(1,0);\draw[line width=0.5,gray] (-2,0)--(2,0);
\draw[line width=1.0] (0,2)--(0,5);\draw[line width=1.0] (1,sqrt 3)--(1,5);\draw[line width=1.0] (1,sqrt 3) arc (60:90:2);
\coordinate[label=center:$T_0$] (o) at (0.5,3);\coordinate[label=left:$i$] (o) at (0,2.2);\coordinate[label=right:${\rho}$] (o) at (1,1.9);
\draw[line width=0.05mm,fill=white](0,2) circle(0.07);\draw[line width=0.05mm,fill=white](1,sqrt 3) circle(0.07);
\end{tikzpicture}&\qquad
\begin{tikzpicture}
[scale=0.6];
\draw[line width=0.5,gray]  (0,0) circle (2);\draw[line width=0.5,gray] (0,5)--(0,-2.5);\draw[line width=1.0] (0,5)--(0,-2);\draw[line width=0.5,gray] (1,5)--(1,-2.5);
\draw[line width=1.0] (1,5)--(1,-sqrt 3);\draw[line width=1.0] (0,-2) arc (270:300:2);
\coordinate[label=center:$T_1$] (o) at (0.5,2);\coordinate[label=left:$-i$] (o) at (0,-2.2);\coordinate[label=right:${\bar\rho}$] (o) at (1,-1.9);
\draw[line width=0.05mm,fill=white](0,-2) circle(0.07);\draw[line width=0.05mm,fill=white](1,-sqrt 3) circle(0.07);
\end{tikzpicture}
\end{tabular}
\end{center} 
\caption{\label{bild_1} {\small The circular triangles $T_0$ and $T_1$.}}
\end{figure}
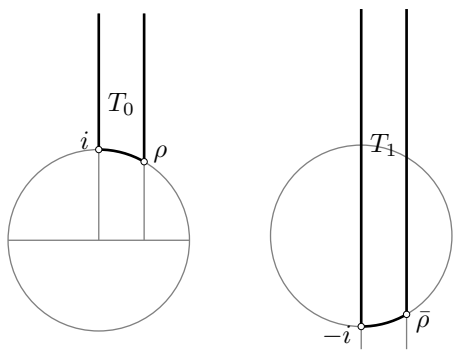

There are essentially two approaches to the {\it elliptic modular function}(\footnote{{\it It is better to solve one problem five different ways, than to solve five different problems one way} (attributed to P\'olya and entirely consistent with his philosophy).}): an analytic approach via the cross-ratio
\be{ellimod1}\lambda(\tau)=(e_1,e_2,e_3,e_4)\quad(\tau=\omega_2/\omega_1)\ee
of the ramified values of the Weierstrass P-function $\wp(z|\omega_1,\omega_2)$ that is associated with the lattice spanned by $\omega_1$ and $\omega_2$,
and a geometric approach via the conformal mapping of the {\it modular triangle} (see Fig.~\ref{bild_5})
\be{modT}\Omega=\{\tau:\Im\tau>0, 0<\Re\tau<1, |\tau|>1\}\ee
onto the upper half-plane $\H$. Both domains are viewed as {\it circular triangles}
with {\it vertices} $0,1,\infty$, so $\lambda:\Omega\longrightarrow\H$ means that the vertices are mapped onto each other in the sense that $f(\tau)\to \xi\in\{0,1,\infty\}$ as $z\to \xi$ in $\Omega$. This is achieved by our non-classical, nevertheless useful definition of the cross-ratio
$$(a,b,c,d)=\frac{(a-c)(b-d)}{(a-d)(b-c)}.$$
In~\cite{steinmetz1} both approaches are  carried out. In \cite{bonk1}, Bonk took the first (classical) path to discuss the ratio
\be{bonkratio}p(\tau)=\eta_1/\eta_2\quad(\tau=\omega_1/\omega_2)\ee
of the quasi-periods $\eta_1$ and $\eta_2$ of the Weierstrass zeta-function
$$\zeta(z)=z^{-1}+\sum_{\omega\ne 0}\big((z-\omega)^{-1}+\omega^{-1}+z\omega^{-2}\big);$$
the sum is over all non-zero lattice points $\omega$, and
$$\zeta(z+\omega_j)=\zeta(z)+\eta_j$$
holds for $j=1,2$. The main result in ~\cite{bonk1}, Theorem 1.1, states that the function (\ref{bonkratio}) maps the domain
$$T_0=\{\tau:\Im\tau>0, 0<\Re\tau<1/2, |\tau|>1\}$$
conformally onto the domain
$$T_1=\{\tau: 0<\Re\tau<1/2,  \Im\tau\ge 0\}\cup\{\tau: \Im\tau<0, |\tau|<1\},$$
see Fig.~\ref{bild_1}; both domains are considered as circular triangles with vertices $i,\rho,\infty$ and $-i,\bar\rho,\infty$ ($\rho=\frac12(1+i\sqrt 3)$)
and interior angles  $\frac12\pi,\frac13\pi,0$ and $\frac12\pi,\frac23\pi,0$, respectively.

\medskip
As explicitly noted in \cite{bonk1}, the results of this expository paper are not new but can be found scattered in the classical literature, sometimes only in implicit form. The paper \cite{bonk1} is written, in part, in the hope that ``it is easily accessible and may serve as an introduction to this classical themes''.
It is easy to agree with this view. To reach a broader audience I've discussed in detail knowledge that experts would take for granted. It is, however, important to note that, in particular in~\cite{bonk2}, the intention was much more ambitious, as the title already suggests, namely to draw conclusions about the critical points of the Eisenstein series $E_2, E_4,$ and $E_6$ via conformal mappings of special circular triangles.
 \medskip

\begin{figure}[!ht]
\begin{center}
\begin{tabular}{ccc}
\begin{tikzpicture}
[scale=0.6];
\draw[line width=0.5,white] (0,0)--(0,-4);\draw[line width=0.5,gray]  (0,0) circle (2);\draw[line width=0.5,gray] (0,5)--(0,0);\draw[line width=0.5,gray] (1,5)--(1,0);\draw[line width=0.5,gray] (-2,0)--(2,0);
\draw[line width=1.0] (0,2)--(0,5);\draw[line width=1.0] (1,sqrt 3)--(1,5);\draw[line width=1.0] (1,sqrt 3) arc (60:90:2);
\coordinate[label=center:$T_0$] (o) at (0.5,3);\coordinate[label=left:$i$] (o) at (0,2.2);\coordinate[label=right:${\rho}$] (o) at (1,1.9);
\draw[line width=0.05mm,fill=white](0,2) circle(0.07);\draw[line width=0.05mm,fill=white](1,sqrt 3) circle(0.07);
\end{tikzpicture}&\qquad
\begin{tikzpicture}
[scale=0.6];
\draw[line width=0.5,gray]  (0,0) circle (2);\draw[line width=0.5,gray] (0,4)--(0,-4);\draw[line width=0.5,gray] (1,4)--(1,-4);
\draw[line width=1.0] (1,4)--(1,sqrt 3);\draw[line width=1.0] (0,-2)--(0,-4);\draw[line width=1.0] (1,sqrt 3) arc (60:270:2);
\coordinate[label=center:$X_0$] (o) at (-2.8,0);\coordinate[label=left:$-i$] (o) at (0,-2.2);\coordinate[label=right:${\rho}$] (o) at (1,1.9);
\draw[line width=0.05mm,fill=white](0,-2) circle(0.07);\draw[line width=0.05mm,fill=white](1,sqrt 3) circle(0.07);
\end{tikzpicture}&\qquad
\begin{tikzpicture}
[scale=0.6];
\draw[line width=0.5,gray]  (0,0) circle (2);\draw[line width=0.5,gray] (0,4)--(0,-4);\draw[line width=0.5,gray] (1,4)--(1,-4);
\draw[line width=1.0] (1,-4)--(1,-sqrt 3);\draw[line width=1.0] (0,4)--(0,2);\draw[line width=1.0] (1,-sqrt 3) arc (-60:90:2);
\coordinate[label=center:$Y_0$] (o) at (2.8,0);\coordinate[label=left:$i$] (o) at (0,2.2);\coordinate[label=right:${\bar\rho}$] (o) at (1,-1.9);
\draw[line width=0.05mm,fill=white](0,2) circle(0.07);\draw[line width=0.05mm,fill=white](1,-sqrt 3) circle(0.07);
\end{tikzpicture}
\end{tabular}
\end{center} 
\caption{\label{bild_2} {\small The circular triangles $T_0$, $X_0$, and $Y_0$.}}
\end{figure}
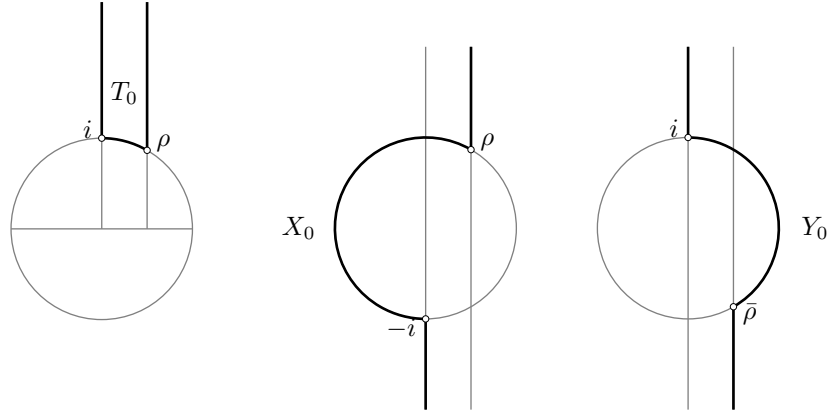

Whilst the Schwarzian derivative takes centre stage here, the approach in \cite{bonk1} uses such classic concepts as the hypergeometric differential equation, series
$$E_2(\tau)=1-24\sum_{k=1}^\infty\frac{n q^k}{1-q^k}\quad(q=e^{2\pi i\tau}, \Im\tau>0)$$
and, later on in \cite{bonk2}, analogues
$$\ds E_{2n}(\tau)=\frac{G_{2n}(\tau)}{2\zeta(2n)}\quad (n\ge 2),$$
where $\zeta$ now denotes the {\it Riemann zeta-function} and the series
$$G_{2n}(\tau)=\sum_{(k,m)\ne (0,0)}(k+m\tau)^{-2n}\quad(n\ge 2, \Im\tau>0)$$
are the classical ones, named after {\it Eisenstein.}(\footnote{The functions $E_{2n}$ involve series $\ds\sum_{k=1}^\infty\frac{k^{2n-1} q^k}{1-q^k}$ with $q=e^{2\pi i\tau},$ $\Im\tau>0$.
Series of the form $\ds\sum_{k=1}^\infty\frac{a_k q^k}{1-q^k}$ are known as {\it Lambert series}.})
These classical functions are used in \cite{bonk2} to find representations of the conformal mappings of the circular triangle $T_0$ onto the circular triangles
$X_0$ and $Y_0$, see~Fig.~\ref{bild_2}. \medskip

Relevant books provide information on all the topics, terms, and facts discussed here: \cite{cara,nehari} on conformal mapping, \cite{bieber} and again \cite{cara} on differential equations,
\cite{hurcour} on elliptic and modular functions, and \cite{steinmetz1} on differential equations, conformal and proper mappings, and elliptic and modular functions.

\section{Circular Triangles and the Schwarzian Derivative}

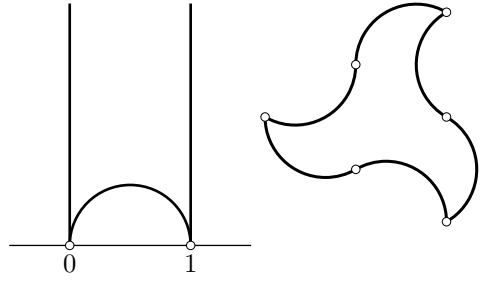
\begin{figure}[!ht]
\begin{center}
\begin{tabular}{cc}
\begin{tikzpicture}
[scale=0.8]
\draw[line width=0.5] (-2,0)--(2,0);\draw[line width=1.0] (-1,0)--(-1,4);\draw[line width=1.0] (1,0)--(1,4);\draw[line width=1.0] (1,0) arc(0:180:1);
\coordinate[label=below:$0$] (o) at (-1,0);\coordinate[label=below:$1$] (o) at (1,0);\draw[line width=0.05mm,fill=white](-1,0) circle(0.07);\draw[line width=0.05mm,fill=white](1,0) circle(0.07);
\end{tikzpicture}
\begin{tikzpicture}
[scale=0.8]
\draw[white] (0,0)--(0,-2.7);
\draw[line width=0.4mm] (1,sqrt 3) arc(60:180:1); \draw[line width=0.4mm] (1,sqrt 3) arc(120:240:1);\draw[line width=0.4mm] (-2,0) arc(180:300:1);
\draw[line width=0.4mm] (1,-sqrt 3) arc(0:120:1);\draw[line width=0.4mm] (1,-sqrt 3) arc(300:420:1);\draw[line width=0.4mm] (-2,0) arc(240:360:1);
\draw[line width=0.05mm,fill=white](1,sqrt 3) circle(0.07);\draw[line width=0.05mm,fill=white](-2,0) circle(0.07);\draw[line width=0.05mm,fill=white](1,-sqrt 3) circle(0.07);
\draw[line width=0.05mm,fill=white](-0.5,-0.5*sqrt 3) circle(0.07);\draw[line width=0.05mm,fill=white](-0.5,0.5*sqrt 3) circle(0.07);\draw[line width=0.05mm,fill=white](1,0) circle(0.07);
\end{tikzpicture}&
\end{tabular}\end{center} 
\caption{\label{bild_5} {\small The modular triangle and the pajarita hexagon (with alternating interior angles $\pi$ and $\frac13\pi$). The latter motif is often found in the Alhambra of Granada.}}
\end{figure}

A circular $n$-gon is a Jordan domain whose boundary curve consists of $n\ge 3$ circular arcs or line segments joining the vertices $z_\nu$ and $z_{\nu+1}$ ($z_{n+1}=z_1$). The boundary curve is oriented (anti-clockwise) in such a way that the domain lies to the left of it. Well-known and lesser-known examples are the {\it modular triangle} and the {\it pajarita hexagon}, see Fig.~\ref{bild_5}.
To any conformal mapping $f$ of the unit disc $\D$ onto some $n$-gon $D$ with interior angles $\pi\alpha_\nu$ at its vertices $z_\nu$ there exists
points $\zeta_\nu$ on $\partial\D$ and constants $C_\nu$ such that the Schwarzian derivative(\footnote{Classically written as $\{f,z\}$, as always with its pros and cons.}) of $f$,
$$\ts S_f=(f''/f')'-\frac12(f''/f')^2$$
has the form
$$S_f(\zeta)=\sum_{\nu=1}^n\Big(\frac{\frac12(1-\alpha_\nu^2)}{(\zeta-\zeta_\nu)^2}+\frac{C_\nu}{\zeta-\zeta_\nu}\Big)$$
with side conditions
$$\sum_{\nu=1}^n C_\nu=\sum_{\nu=1}^n(\zeta_\nu C_\nu+{\ts\frac12(1-\alpha_\nu^2)})=\sum_{\nu=1}^n(\zeta^2_\nu C_\nu+(1-\alpha_\nu^2)\zeta_\nu)=0.$$
Three of the `vertices' $\zeta_\nu$ may be prescribed, and only in particularly symmetric situations(\footnote{The pajarita hexagon is not enough symmetric: the Schwarzian derivative on the unit disc $\D$ has the form $\ds\frac{4z(az^3+b)}{(1-z^3)^2(c^3-z^3)}$ with $c^3-1=a+b$; $a$ and $b$ have to be computed numerically.}) or in case of a {\it circular triangle} ($n=3$) the constants $C_\nu$
and the points $\zeta_\nu$ ($4\le\nu\le n$) can be determined from the side conditions (with Vandermonde determinant) in terms of the $\alpha_\nu$'s.
In the usual sense, any two circular triangles $T$ and $T'$ are conformally equivalent, that is, there always exists some conformal mapping $h:T'\kzu1{1}T$. This is even true if we demand that vertices are mapped onto vertices.
In particular, the upper half-plane $\H$ with vertices $0,1,\infty$ can be mapped conformally onto the circular triangle $T$ such that $0,1,\infty$ correspond to the vertices of $T$ in anti-clockwise order. With the help of the already mentioned side conditions we obtain that

\begin{figure}[!ht]
\begin{center}
\begin{tabular}{cc}
\begin{tikzpicture}
[scale=0.6];
\draw[line width=0.5,gray]  (0,0) circle (2);\draw[line width=0.5,gray] (-2,-2)--(-2,4);\draw[line width=0.5,gray] (-1,-2)--(-1,4);
\draw[line width=1.0] (-2,0)--(-2,4);\draw[line width=1.0] (-1,sqrt 3)--(-1,4);\draw[line width=1.0] (-1,sqrt 3) arc (120:180:2);
\coordinate[label=center:$U_0$] (o) at (-1.5,3);\coordinate[label=left:$0$] (o) at (-2,0);\coordinate[label=right:${\rho}$] (o) at (-1,2.1);
\draw[line width=0.05mm,fill=white](-2,0) circle(0.07);\draw[line width=0.05mm,fill=white](-1,sqrt 3) circle(0.07);
\end{tikzpicture}&\qquad
\begin{tikzpicture}
[scale=0.6];
\draw[line width=0.5,gray]  (0,0) circle (2);\draw[line width=0.5,gray] (-2,-2)--(-2,4);\draw[line width=0.5,gray] (-1,-2)--(-1,4);
\draw[line width=1.0] (-2,0)--(-2,4);\draw[line width=1.0] (-1,-sqrt 3)--(-1,4);\draw[line width=1.0] (-2,0) arc (180:240:2);
\coordinate[label=left:$0$] (o) at (-2,0);\coordinate[label=left:${\bar\rho}$] (o) at (-1,-1.9);\coordinate[label=center:$Z_0$] (o) at (-1.5,1.5);
\draw[line width=0.05mm,fill=white] (-2,0) circle(0.07);\draw[line width=0.05mm,fill=white](-1,-sqrt 3) circle(0.07);
\end{tikzpicture}
\end{tabular}
\end{center} 
\caption{\label{bild_3} {\small The circular triangles $U_0$ and $Z_0$.}}
\end{figure}
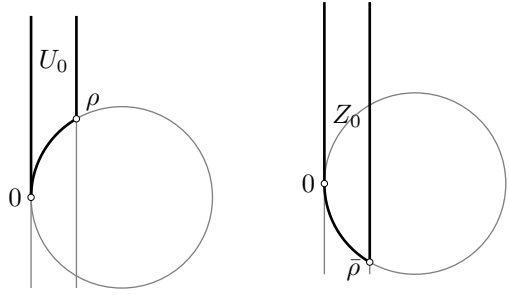

\medskip\begin{itemize}\item[$\vartriangleright$] {\it the conformal mapping $f$ that maps the upper half-plane $\H$ onto some circular triangle $T$ and the points $0,1,\infty$ onto the 
vertices $z_1,z_2,z_3$ of $T$ has the Schwarzian derivative{\rm(\footnote{This is, of course, classical. It takes three minutes to write a
\texttt{maple} code for computing $S_f$ in terms of the $\alpha_\nu$'s. Life is so much easier today, but at the same time we are losing skills.})}
\be{S3alpha}S_f(z)=\frac{\frac12(1-\alpha_1^2)}{z^2}+\frac{\frac12(1-\alpha_2^2)}{(z-1)^2}+\frac{\frac12(\alpha_1^2+\alpha_2^2-\alpha_3^2-1)}{z(z-1)}.\ee%
In particular, the Schwarzian derivatives of the conformal mappings of $\H$ onto various circular triangles $T_0, T_1, U_0, X_0, Y_0, Z_0$ occurring in \cite{bonk1,bonk2}
and also here {\rm (see Fig.~\ref{bild_1}, \ref{bild_2} and \ref{bild_3})} are given by}

\begin{tabular}{crl}
$(T_0)$& $\ds\frac{\frac49}{z^2}+\frac{\frac38}{(z-1)^2}-\frac{\frac{23}{72}}{z(z-1)}=$&$\!\!\!\ds\frac{36z^2-31z+27}{72z^2(z-1)^2}$\cr
$(T_1)$& $\ds\frac{\frac5{18}}{z^2}+\frac{\frac38}{(z-1)^2}-\frac{\frac{11}{72}}{z(z-1)}=$&$\!\!\!\ds\frac{36z^2+43z+27}{72z^2(z-1)^2}$\cr
$(X_0)$& $\ds\frac{\frac49}{z^2}+\frac{\frac38}{(z-1)^2}-\frac{\frac{59}{72}}{z(z-1)}=$&$\!\!\!\ds\frac{5z+27}{72z^2(z-1)^2}$\cr
$(Y_0)$& $\ds\frac{\frac5{18}}{z^2}+\frac{\frac38}{(z-1)^2}-\frac{\frac{47}{72}}{z(z-1)}=$&$\!\!\!\ds\frac{-7z+27}{72z^2(z-1)^2}$\cr
\end{tabular}

\begin{tabular}{crl}
$(U_0)$& $\ds\frac{\frac12}{z^2}+\frac{\frac5{18}}{(z-1)^2}-\frac{\frac{5}{18}}{z(z-1)}=$&$\!\!\!\ds\frac{9z^2-13z+9}{72z^2(z-1)^2}$\cr
$(Z_0)$& $\ds\frac{\frac49}{(z-1)^2}+\frac{\frac1{18}}{z(z-1)}=$&$\!\!\!\ds\frac{9z-1}{18z(z-1)^2}$
\end{tabular}
\end{itemize}


The Schwarzian derivative is invariant under the full M\"obius group, that is,
$$S_{M\circ f}=S_f$$
holds for every M\"obius transformation $\ds M=\frac{az+b}{cz+d}$ ($ad-bc\ne 0$), and $S_f=0$ if and only if $f$ is a M\"obius transformation. More generally we have

\medskip\begin{itemize}
\item[$\vartriangleright$] {\it the chain rule\quad $S_{g\circ f}=(S_g\circ f)f'^2+S_f.$}
\end{itemize}

In particular, if $f$ and $g$ are local inverse function of each other, the chain rule yields $0=(S_g\circ f)f'^2+S_f$, hence
\be{CRinv}S_{f}=-(S_{f^{-1}}\circ f)f'^2.\ee
If $f_1$ and $f_2$ map the circular triangle $T_1$ onto $\H$ and $\H$ onto the circular triangle $T_2$, respectively, then $h=f_2\circ f_1$ maps $T_1$ onto $T_2$ with Schwarzian derivative
\be{hT1T2}S_h=(S_{f_2}\circ f_1-S_{f^{-1}_1}\circ f_1)f_1'^2.\ee
The highly non-linear differential equation
\be{SDGL}S_f(z)=2Q(z)\ee
can be reduced to the second order linear differential equation
\be{LDGL}w''+Q(z)w=0\ee
in the following sense:

\medskip\begin{itemize}
\item[$\vartriangleright$] {\it Given any (local) solution to $(\ref{SDGL})$ there exist linearly independent solutions $w_1$ and $w_2$ to $(\ref{LDGL})$ such that $f=w_1/w_2$ holds, and,
{\rm vice versa}, for any two linearly independent solutions we have $S_{w_1/w_2}(z)=2Q(z)$.}
\end{itemize}

In case of $n$-gon mappings, this leads to differential equations (\ref{LDGL}) of Fuchsian type. For $n=3$, in particular, (\ref{LDGL}) may be transformed into the
{\it hypergeometric differential equation}
\be{HDGL}z(1-z)w''+(\gamma-(\alpha+\beta+1)z)w'+\alpha\beta w=0,\ee
with parameters $\alpha,\beta,\gamma$ corresponding to the numbers $\alpha_1,\alpha_2,\alpha_3$.

If $T$ is any circular triangle with vertices $z_\nu$ and interior angles $\pi\alpha_\nu$, then  any M\"obius transformation $M$  maps $T$ conformally onto
the circular triangle $T^*$ with vertices $z^*_\nu=M(z_\nu)$ and interior angles $\pi\alpha_\nu$;
$T$ and $T^*=M(T)$ are called M\"obius equivalent. Now suppose $T$ and $T^*$ are circular triangles with vertices $z_\nu,z^*_\nu$ and interior angles $\pi\alpha_\nu=\pi\alpha1*_\nu$, and
let $f$ and $h$ denote the conformal mappings of $\H$ onto $T$ and $T$ onto $T^*$ that send $0,1,\infty$ to $z_1,z_2,z_3$ and $z_1,z_2,z_3$ to $z^*_1,z^*_2,z^*_3$ (both in this order)
respectively. Then $f$ and $h\circ f$ have the same Schwarzian derivative, hence
$$S_f=S_{h\circ f}=(S_h\circ f)f'^2+S_f,$$
that is, $S_h=0$ holds on $T$ and $h$ is a M\"obius transformation. Thus

\medskip\begin{itemize}\item[$\vartriangleright$] {\it any two circular triangles are M\"obius equivalent if and only if they have the same interior angles.}\end{itemize}
\section{The Absolute Invariant}
The Weierstrass P-function $\wp$ satisfies
$$\wp'^2=4\wp^3-g_2\wp-g_3=4(\wp-e_1)(\wp-e_2)(\wp-e_3)$$
with $g_2=60G_4$ and $g_3=140G_6$ and non-vanishing discriminant $g_2^3-27g_3^2$. The absolute invariant or $j$-function
$$\j(\tau)=\frac{g_2(\tau)^3}{g_2(\tau)^3-27g_3(\tau)^2}$$
(as in~\cite{steinmetz1}, with some advantage, the traditional factor $12^3=1728$ is omitted) is holomorphic on the upper half-plane, invariant under the full modular group PSL$_2(\Z)$ and a rational function of degree $6$ of the elliptic modular function.
\begin{figure}[!ht]
\begin{center}\begin{tabular}{cc}
\begin{tikzpicture}
[scale=0.4];
\coordinate[label=below:$-1$] (-1) at (0,0);\coordinate[label=below:$0$] (0) at (4,0);\coordinate[label=below:$1$] (1) at (8,0);
\draw[line width=0.1mm,fill=gray!50!] (0,0) arc(180:90:2) -- (2,3.464) arc(120:180:4);\draw[line width=0.1mm,fill=gray!50!] (4,0) arc(180:90:2) -- (6,3.464) arc(120:180:4);
\draw[line width=0.1mm,fill=gray!20!] (2,2) arc(90:0:2) arc(0:60:4) --(2,2);\draw[line width=0.1mm,fill=gray!20!] (6,2) arc(90:0:2) arc(0:60:4) --(6,2);
\draw[line width=0.1mm,fill=gray!20!] (0,0) arc(180:120:4) arc(60:90:4) --(0,0);\draw[line width=0.1mm,fill=gray!20!] (4,0) arc(180:120:4) arc(60:90:4) --(4,0);
\draw[line width=0.1mm,fill=gray!50!] (4,0)--(4,4) arc(90:120:4) arc(60:0:4);\draw[line width=0.1mm,fill=gray!50!] (8,0)--(8,4) arc(90:120:4) arc(60:0:4);
\draw[line width=0.1mm,fill=gray!50!] (2,8)--(2,3.464) arc(60:90:4) --(0,8);\draw[line width=0.1mm,fill=gray!50!] (6,8)--(6,3.464) arc(60:90:4) --(4,8);
\draw[line width=0.1mm,fill=gray!20!] (2,8)--(2,3.464) arc(120:90:4) --(4,8);\draw[line width=0.1mm,fill=gray!20!] (6,8)--(6,3.464) arc(120:90:4) --(8,8);
\draw[line width=0.4mm,gray] (6,8)--(6,2*sqrt 3);\draw[line width=0.4mm,gray] (6,2*sqrt 3) arc (120:180:4);\draw[line width=0.4mm,gray] (8,0) arc (0:60:4);
\draw[line width=0.4mm,gray] (2,8)--(2,2*sqrt 3);\draw[line width=0.4mm,gray] (2,2*sqrt 3) arc (120:180:4);\draw[line width=0.4mm,gray] (4,0) arc (0:60:4);
\draw[line width=0.4mm,white] (8,4) arc(90:120:4);\draw[line width=0.4mm,white] (4,4) arc(90:120:4);
\draw[line width=0.4mm,white] (6,2)--(6,2*sqrt 3);\draw[line width=0.4mm,white] (2,2)--(2,2*sqrt 3);
\draw[line width=0.4mm,white] (6,2*sqrt 3) arc(60:90:4);\draw[line width=0.4mm,white] (2,2*sqrt 3) arc(60:90:4);
\draw[line width=0.4mm,black] (8,0)--(8,8);\draw[line width=0.4mm,black] (4,8)--(4,0);\draw[line width=0.4mm,black] (8,0) arc(0:180:2);\draw[line width=0.4mm,black] (4,0) arc(0:180:2);
\coordinate[label=center:$T_0$] (omega) at (5,6);\coordinate[label=center:$T'_0$] (omega) at (4.6,3.2);
\draw[line width=0.4mm] (0,0)--(0,8);\draw[line width=0.4mm] (4,0)--(4,8);\draw[line width=0.4mm] (8,0)--(8,8);
\draw[line width=0.2mm,fill=white] (0,0) circle (0.125);\draw[line width=0.2mm,fill=white] (4,0) circle (0.125);\draw[line width=0.2mm,fill=white] (8,0) circle (0.125);
\draw[line width=0.2mm,fill=white] (2,2*sqrt 3) circle (0.125);\draw[line width=0.2mm,fill=white] (6,2*sqrt 3) circle (0.125);
\draw[line width=0.2mm,fill=white] (0,4) circle (0.125);\draw[line width=0.2mm,fill=white] (4,4) circle (0.125); \draw[line width=0.2mm,fill=white] (8,4) circle (0.125);
\draw[line width=0.2mm,fill=white] (2,2) circle (0.125);\draw[line width=0.2mm,fill=white] (6,2) circle (0.125);
\end{tikzpicture}&\qquad\qquad
\begin{tikzpicture}
[scale=0.55];
\draw[white] (0,0)--(0,-2.8);
\draw[line width=0.1mm,fill=gray!20!,gray!20!] (-2,0) rectangle (5,2);\draw[line width=0.4mm,fill=gray!50!,gray!50!] (-2,-2) rectangle (5,0);
\draw[line width=0.4mm,black] (5,0.05)--(2,0.05);
\draw[line width=0.4mm,white] (2,0.05)--(0,0.05);
\draw[line width=0.4mm,gray] (-2,0.05)--(0,0.05);
\coordinate[label=below:$0$] (o) at (0,0);\coordinate[label=below:$1$] (o) at (2,0);
\draw[line width=0.2mm,fill=white] (0,0) circle (0.085);\draw[line width=0.2mm,fill=white] (2,0) circle (0.085);
\end{tikzpicture}
\end{tabular}
\end{center} 
\caption{\label{bild_4} {\small The circular arcs and vertical lines bound $12$ circular triangles, each of which is mapped by the $\j$-function alternately onto the upper and lower half-plane, respectively.}
}\end{figure}
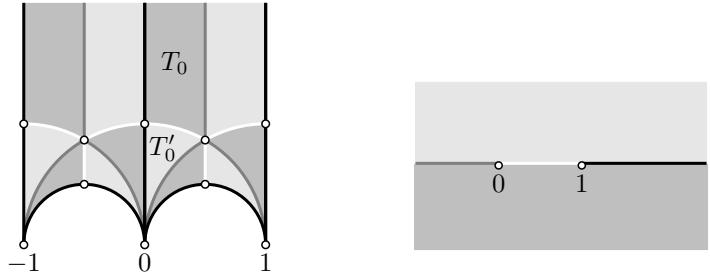
 This rational function
$$\ds R(\zeta)=\frac 4{27}\frac{(1-\zeta+\zeta^2)^3}{\zeta^2(1-\zeta)^2}$$
is invariant under the {\it cross-ratio group}, consisting of those M\"obius transformations that permute the points $0,1,\infty$, and a ramified proper mapping(\footnote{A proper mapping $f:D\longrightarrow G$ between planar domains $D$ and $G$ is a holomorphic function that maps $\partial D$ into $\partial G$ in the sense of $f(z)\to\partial G$ as $z\to\partial D$. Proper mappings have a mapping degree $d$, that is, each $w\in G$ has $d$ pre-images (counted by multiplicities) in $D$. Conformal mappings are proper mappings of degree $1$, and {\it vice versa}.}) of degree $6$ of $\C\setminus\{0,1\}$ onto $\C\setminus[1,+\infty)$. The $\j$-function itself is a proper mapping of the modular triangle~(\ref{modT}) of degree $3$ onto $\C\setminus[1,+\infty)$ and maps the circular triangle $T_0$ conformally onto the lower half-plane $\H_-$ and the vertices $i,\rho,\infty$ onto $1,0,\infty$ (see Fig.~\ref{bild_4}). It is therefore appropriate, as already noted in \cite{bonk1,bonk2} to consider
$h:T_0\longrightarrow D$, where $D=T_1$ or $X_0$ or $Y_0$ (or any other circular triangle) as a function of $\j$, that is, consider $f=h\circ\j^{-1}:\H_-\longrightarrow D$.
Then (\ref{hT1T2}) with $f_1=\j$ and $f_2=f$ and yields
$$S_h=([S_{f}-S_{\j^{-1}}]\circ\j)\j'^2,$$
and in more detail, using the notation from \cite{bonk1,bonk2}
$$\begin{aligned}
p(\tau)=&~\tau-\frac{6i}{\pi E_2(\tau)}\cr
s_4(\tau)=&~\tau-\frac{6i E_4(\tau)}{\pi(E_2(\tau)E_4(\tau)-E_6(\tau))}\cr
s_6(\tau)=&~\tau-\frac{6i E_6(\tau)}{\pi(E_2(\tau)E_6(\tau)-E_4^2(\tau))}
\end{aligned}$$

\medskip\begin{itemize}
\item[$\vartriangleright$] {\it The  Schwarzian derivatives of the conformal mappings
$$p:T_0\longrightarrow T_1,~s_4:T_0\longrightarrow X_0,~{\rm and~}s_6:T_0\longrightarrow Y_0$$
have the form 
$$S_p=\frac{\j'^2}{6\j^2(\j-1)},~ S_{s_4}=\frac{-\j'^2}{2\j(\j-1)}, ~{\rm and}~ S_{s_6}=\frac{(1-3\j)\j'^2}{6\j^2(\j-1)},$$
respectively. Each mapping admits unrestricted analytic continuation in the upper half-plane $\H$
\begin{itemize}
\item[-] either by repeated application of the Schwarz reflection principle
\item[-] or else as quotients of {\rm (globally existing!)} linearly independent solutions to the respective linear differential equation
$$w''+\ts\frac12 S_h(\tau)w=0$$
\end{itemize}
$(h=p, s_4, s_6)$; in each case
$$h\circ M=M\circ h$$
holds for every M\"obius transformation $\ds M(\tau)=\frac{a\tau+b}{c\tau+d}$ in {\rm PSL}$_2(\Z)$. In particular, $S_h$ is a modular form of weight $4$,
$$S_h(M(\tau))=(c\tau+d)^{4}S_h(\tau).$$}
\end{itemize}

Since the group PSL$_2(\Z)$ is generated by the transformations $\tau\mapsto\tau+1$ and $\tau\mapsto -1/\tau$, the last assertion follows from the reflection principle:
$$h(-\bar \tau)=-\overline{h(\tau)},~h(1-\bar \tau)=1-\overline{h(\tau)},~{\rm and~}h(1/\bar \tau)=1/\overline {h(\tau)},$$
thus  $h(1-\bar \tau)=1+h(-\bar \tau)$, that is $h(\tau+1)=h(\tau)+1$ (set $\tau=-\bar\tau$) and, in the same manner, $h(-1/\tau)=-1/h(\tau)$. This holds throughout $\H$ by the permanence principle.
In any case, $h(\tau)-\tau$ is $1$-periodic,
$$h(\tau)=\tau+\sum_{n\in\Z}c_ne^{2\pi i \tau}$$
holds on the half-plane $\Im\tau>1$, where all these functions are pole-free. Since $\Im h(\tau)\to +\infty$ as $\Im\tau\to\infty$, the coefficients $c_n$ with $n<0$ vanish
and
$$h(\tau)=\tau+\sum_{n=0}^\infty c_ne^{2\pi i \tau}\quad(\Im\tau>1)$$
holds. Since $p$ doubles angles at $\tau\equiv \rho$ modulo PSL$_2(\Z)$, $p'$ vanishes at these points, and nowhere else.
We note that $\eta_1$ and $\eta_2$ do {\it not} solve $w''+\frac12 S_p(\tau)w=0$. Instead, the following applies: $p=w_1/w_2$ with $w_j=\eta_j/W$, where $W$ denotes the Wronskian
determinant 
$$\left|\!\!\begin{array}{cc}\eta_1&\eta_2\cr\eta_1'&\eta_2'\end{array}\!\!\right|=\eta_2^2-2\pi i\eta_2',$$ 
and $w_1, w_2$ form a fundamental system of $w''+\frac12 S_p(\tau)w=0$.

\section{One More Example}

The $\j$-function maps $U_0$ conformally onto $D=\C\setminus((-\infty,0]\cup[1,+\infty))$, see Fig.~\ref{bild_4}, with Schwarzian derivative $S_{\j^{-1}}(z)=\ds\frac{36z^2-31z+27}{72z^2(z-1)^2}$, $z=k(u)=1/(1-u^2)$ maps the upper half-plane $\H$ conformally onto $D$ with $S_k(u)=\ds -\frac{\frac32}{u^2}$, and $f$ maps $\H$ onto $Z_0$ with Schwarzian derivative $S_f(u)=\ds\frac{9u-1}{18u(u-1)^2}$.
Set $\phi=k^{-1}\circ \j:U_0\longrightarrow \H$. Then $h=f\circ \phi$ maps $U_0$ conformally onto $Z_0$ and
$$\begin{aligned}
\phi'^2=&~\frac{\j'^2}{4(\j-1)\j^3}\cr
S_h=&~(S_f\circ\phi)\phi'^2+S_\phi\cr
S_\phi=&~-(S_{\phi^{-1}}\circ\phi)\phi'^2\cr
S_{\phi^{-1}}=&~(S_{\j^{-1}}\circ k) k'^{2}+S_k
\end{aligned}$$
hold, and from
$$S_h=([S_f-(S_{\j^{-1}}\circ k)k'^2-S_k]\circ \phi)\phi'^2$$
it follows that (with the notation in \cite{bonk2})\medskip

\begin{itemize}\item[$\vartriangleright$] \quad $\ds S_{s_2^\pm}=\frac{-20\j^2+6\j+9\pm\sqrt{\j-1}(6\j-17)}{72\j^2(\j-1)^2}\j'^2$;
\end{itemize}
there are two holomorphic branches of $\sqrt{\j-1}$ on $\H$, see also (10.8) for $s_2^+$ and (10.9) for $s_2^-$ in~\cite{bonk2}.

\bigskip{\footnotesize Institut f\"ur Mathematik\\ Technische Universit\"at Dortmund\\
stein@math.tu-dortmund.de\\ www.mathematik.tu-dortmund.de/steinmetz/\\
www.researchgate.net/profile/Norbert-Steinmetz}
\end{document}